\documentclass{article}

\usepackage{graphicx}

\usepackage{amsfonts}

\usepackage{hyperref}

\usepackage{mathtools}

\usepackage{cleveref}

\usepackage{tikz}
\usepackage{xcolor}

\usepackage[left=25mm, right=25mm, bottom=35mm]{geometry}

\newcommand{\Params}{\mathcal{P}}
\newcommand{\setR}{\mathbb{R}}
\newcommand{\setN}{\mathbb{N}}
\newcommand{\dt}{\,\mathrm{d}t}
\newcommand{\ds}{\,\mathrm{d}s}
\newcommand{\Gramian}{\Lambda_\mu^R}
\DeclareMathOperator*{\argmin}{arg\,min}

\newcommand{\norm}[1]{\left\lVert#1\right\rVert}

\definecolor{color0}{rgb}{0.65,0,0.15}
\definecolor{color1}{rgb}{0.84,0.19,0.15}
\definecolor{color2}{rgb}{0.96,0.43,0.26}
\definecolor{color3}{rgb}{0.99,0.68,0.38}
\definecolor{color4}{rgb}{1,0.88,0.56}
\definecolor{color5}{rgb}{0.67,0.85,0.91}
\definecolor{color6}{rgb}{0.27,0.46,0.71}

\usetikzlibrary{shapes.geometric, arrows, positioning, calc, decorations.pathreplacing, tikzmark, patterns}
\pgfdeclarelayer{background}
\pgfdeclarelayer{foreground}
\pgfsetlayers{background,main,foreground}

\tikzstyle{input} = [circle, minimum height=.6cm, text centered, draw=color3, fill=color3!30, font=\small]
\tikzstyle{model} = [rounded corners=0.35cm, minimum height=.7cm, text centered, draw=color1, fill=color1!30, font=\small]
\tikzstyle{error-estimator} = [diamond, text centered, draw=color5, fill=color5!30, aspect=2, inner sep=-.5pt, font=\small]
\tikzstyle{output} = [rectangle, minimum width=2.6cm, minimum height=.7cm, text centered, draw=color6, fill=color6!30, font=\small]
\tikzstyle{arrow} = [thick, ->, >=stealth]

\newcommand{\distbg}{0.15}

\newcommand{\background}[5]{%
	\begin{pgfonlayer}{background}
		\path (#1.west |- #2.north)+(-\distbg,\distbg) node (a1) {};
		\path (#3.east |- #4.south)+(\distbg,-\distbg) node (a2) {};
		\path[fill=color4!30, rounded corners, draw=color4, dashed, thick] (a1) rectangle (a2);
		\path (#3.east |- #2.north)+(0,\distbg)--(#1.west |- #2.north)+(0,\distbg) node[midway] (#5-n) {};
		\path (#3.east |- #2.south)+(0,-\distbg)--(#1.west |- #2.south) node[midway] (#5-s) {};
		\path (#1.west |- #2.north)+(\distbg,0)--(#1.west |- #4.south) node[midway] (#5-w) {};
		\path (#3.east |- #2.north)+(\distbg,0)--(#3.east |- #4.south) node[midway] (#5-e) {};
		
		\node (#5-ne) at ($(#3.east |- #2.north)+(\distbg,\distbg)$) {};
		\node (#5-se) at ($(#3.east |- #4.south)+(\distbg,-\distbg)$) {};
		\node (#5-nw) at ($(#1.west |- #2.north)+(-\distbg,\distbg)$) {};
		\node (#5-sw) at ($(#1.west |- #4.south)+(-\distbg,-\distbg)$) {};
	\end{pgfonlayer}%
}

\newcommand{\FOM}{FOM}
\newcommand{\RBROM}{RB-ROM}
\newcommand{\MLROM}{ML-ROM}

\usepackage{booktabs}

\usepackage{dcolumn}
\newcolumntype{d}[1]{D{.}{.}{#1}}
\newcommand\mc[1]{\multicolumn{1}{c}{#1}}

\usepackage{makecell}

\usepackage{numprint}
\npthousandsep{,}

\usepackage{orcidlink}

\title{Application of an adaptive model hierarchy\\to parametrized optimal control problems\thanks{Funded by the Deutsche Forschungsgemeinschaft (DFG, German Research Foundation) under Germany's Excellence Strategy EXC 2044 –390685587, Mathematics Münster: Dynamics–Geometry–Structure.}}

\author{Hendrik Kleikamp\,\orcidlink{0000-0003-1264-5941}\thanks{Institute for Analysis and Numerics, Mathematics Münster, University of Münster, Einsteinstrasse 62, 48149 Münster, Germany ({\tt hendrik.kleikamp@uni-muenster.de}).}}

\begin{document}

\maketitle

\begin{abstract}
\noindent In this contribution we apply an adaptive model hierarchy, consisting of a full-order model, a reduced basis reduced order model, and a machine learning surrogate, to parametrized linear-quadratic optimal control problems. The involved reduced order models are constructed adaptively and are called in such a way that the model hierarchy returns an approximate solution of given accuracy for every parameter value. At the same time, the fastest model of the hierarchy is evaluated whenever possible and slower models are only queried if the faster ones are not sufficiently accurate. The performance of the model hierarchy is studied for a parametrized heat equation example with boundary value control.
\end{abstract}

\pagestyle{myheadings}
\thispagestyle{plain}
\markboth{H. KLEIKAMP}{Adaptive model hierarchy for parametrized optimal control problems}

\section{Introduction}
Optimal control problems with parameter-dependent system components typically require an enormous computational effort when considered in a multi-query or real time scenario. Solving these kinds of problems exactly for many different values of the parameter is computationally demanding and often prohibitively costly. In~\cite{lazar2016greedy}, a greedy procedure to construct a reduced order model for parametrized optimal control problems has been proposed. To further speed up the online computations of the reduced model, machine learning algorithms have been used in a certified manner in~\cite{kleikamp2023greedy}. In this work, we combine the aforementioned ideas with an adaptive and certified model hierarchy for parametrized problems which was introduced in~\cite{haasdonk2023certified} and further applied in~\cite{wenzel2023application}. This model hierarchy allows for an adaptive construction and improvement of reduced order models and machine learning surrogates while already querying the model hierarchy for different parameter values. Hence, no costly offline phase is required, while the results provided by the model hierarchy still fulfill a prescribed error tolerance. The parametrized optimal control problem is only solved exactly using the underlying full-order model if necessary. Whenever possible, cheaper reduced order models are used, which are built, trained and improved on the fly using data from more accurate but at the same time more costly models.
\par
The adaptive model hierarchy for optimal control problems introduced in this contribution could for instance be applied in conjunction with Monte Carlo estimation of derived quantities. In addition, parameter optimization problems with optimal control problems as constraint could be another possible field of usage.
\par
Applications of reduced order models to optimal control problems can for instance be found in~\cite{ballarin2022spacetime,dede2012reduced,lazar2016greedy}. A combination with tools from machine learning has been proposed in~\cite{daniel2020model}. In~\cite{hesthaven2018nonintrusive}, an approach using deep neural networks for the solution of parametrized partial differential equations combined with reduced basis reduced order models constructed via proper orthogonal decomposition has been discussed.
\par
The paper is organized as follows: In~\Cref{sec:optimal-control-problem}, we introduce the problem considered in this contribution and present the associated optimality system. Afterwards, in~\Cref{sec:reduced-order-models}, the two reduced order models and an a posteriori error estimator will be discussed. \Cref{sec:adaptive-model-hierarchy} presents the adaptive model hierarchy in a general formulation which is then applied to the parametrized optimal control setting in~\Cref{sec:application-model-hierarchy-optimal-control}. A numerical example showing the performance of the devised algorithms is performed and evaluated in~\Cref{sec:numerical-experiment}. The paper ends in~\Cref{sec:conclusion-outlook} with some concluding remarks and an outlook to future research directions.

\section{Linear-quadratic parametrized optimal control problems}\label{sec:optimal-control-problem}
First, we introduce the parametrized optimal control problems considered in this work. For simplicity, we state the optimal control problem in a finite-dimensional setting. The more general formulation for infinite-dimensional parameter, state and control spaces can be found in~\cite{kleikamp2023greedy}. Afterwards, the optimality system using an adjoint variable and a linear system of equations for the optimal adjoint at final time are presented.

\subsection{Problem formulation}
Let~$\Params\subset\setR^p$ be a compact parameter set for some~$p\in\setN$. For a parameter~$\mu\in\Params$, the state system is given as
\begin{equation}\label{equ:state-system}
	\begin{aligned}
		\dot{x}_\mu(t) &= A_\mu x_\mu(t) + B_\mu u(t),\qquad t\in[0,T], \\
		x_\mu(0) &= x_\mu^0,
	\end{aligned}
\end{equation}
where~$x_\mu\colon[0,T]\to\setR^n$ is the state trajectory, $u\colon[0,T]\to\setR^m$ denotes the control, $A_\mu\in\setR^{n\times n}$ is the state operator, $B_\mu\in\setR^{n\times m}$ the control operator, $T>0$ is the final time and~$x_\mu^0\in\setR^n$ denotes the initial state. In the examples we have in mind, for instance discretizations of time-dependent partial differential equations (PDEs), the dimension~$n\in\setN$ of the state space is typically quite large whereas the number of controls~$m\in\setN$ is of moderate size. For each parameter~$\mu\in\Params$, we aim to steer the system state~$x_\mu(T)$ at time~$T$ close to a given target state~$x_\mu^T\in\setR^n$ while not spending too much control energy. We hence aim to minimize the following functional~$\mathcal{J}_\mu$ defined for a control~$u\colon[0,T]\to\setR^m$ as
\[
	\mathcal{J}_\mu(u) \coloneqq \frac{1}{2}\left[(x_\mu(T)-x_\mu^T)^\top M(x_\mu(T)-x_\mu^T)+\int\limits_0^T u(t)^\top Ru(t)\dt\right],
\]
where~$x_\mu\colon[0,T]\to\setR^n$ solves the state equation~\eqref{equ:state-system} for the control~$u$, the matrix~$M\in\setR^{n\times n}$ is symmetric and positive-semidefinite and the matrix~$R\in\setR^{m\times m}$ is symmetric and positive-definite. The matrices~$M$ and~$R$ allow for individual weights for different state and control components. To summarize, given a parameter~$\mu\in\Params$, we would like to solve the optimal control problem
\begin{align}\label{equ:optimal-control-problem}
	\min\limits_{u}\mathcal{J}_\mu(u),\quad\text{s.\,t. }\dot{x}_\mu(t)=A_\mu x_\mu(t)+B_\mu u(t)\text{ for }t\in[0,T],\quad x_\mu(0)=x_\mu^0.
\end{align}

\subsection{Optimality system}\label{subsec:optimality-system}
The following optimality system characterizes the optimal state trajectory~$x_\mu^*\colon[0,T]\to\setR^n$, the optimal control~$u_\mu^*\colon[0,T]\to\setR^m$ and the optimal adjoint trajectory~$\varphi_\mu^*\colon[0,T]\to\setR^n$ (see~\cite[Theorem~2.4]{kleikamp2023greedy} for more details) that solve~\eqref{equ:optimal-control-problem}:
\begin{subequations}\label{equ:optimality-system-main}
	\begin{equation}\label{equ:optimality-system-odes}
		\begin{aligned}
			\dot{x}_\mu^*(t) &= A_\mu x_\mu^*(t)+B_\mu u_\mu^*(t), \\
			-\dot{\varphi}_\mu^*(t) &= A_\mu^\top \varphi_\mu^*(t), \\
			u_\mu^*(t) &= -R^{-1}B_\mu^\top \varphi_\mu^*(t),
		\end{aligned}
	\end{equation}
	for~$t\in[0,T]$ with initial respectively terminal conditions
	\begin{align}\label{equ:optimality-system-boundary-conditions}
		x_\mu^*(0) = x_\mu^0,\qquad \varphi_\mu^*(T)=M\left(x_\mu^*(T)-x_\mu^T\right).
	\end{align}
\end{subequations}
Using the exponential function, we define the \emph{weighted controllability Gramian}~$\Gramian\in\setR^{n\times n}$ as
\[
	\Gramian \coloneqq \int\limits_0^T e^{A_\mu(T-s)}B_\mu R^{-1}B_\mu^\top e^{A_\mu^\top(T-s)}\ds.
\]
Using this Gramian matrix, one can derive a linear system for the optimal final time adjoint~$\varphi_\mu^*(T)$ (see~\cite[Lemma~2.5]{kleikamp2023greedy}), which is given as
\begin{align}\label{equ:linear-system}
	\left(I+M\Gramian\right)\varphi_\mu^*(T) = M\left(e^{A_\mu T}x_\mu^0-x_\mu^T\right).
\end{align}
The solution of the optimality system~\eqref{equ:optimality-system-main} is already uniquely determined by the optimal final time adjoint~$\varphi_\mu^*(T)$ (in the following we assume that the product~$M\Gramian$ is positive-semidefinite for all parameters~$\mu\in\Params$). It is therefore sufficient to first solve the linear system in~\eqref{equ:linear-system} for~$\varphi_\mu^*(T)$ and afterwards solve the ordinary differential equation system in~\eqref{equ:optimality-system-odes} to obtain the optimal control, state and adjoint trajectories. The linear system can be solved using iterative methods since applying the Gramian matrix~$\Gramian$ to a vector~$p\in\setR^n$ is (up to a minus sign) equivalent to solving the optimality system~\eqref{equ:optimality-system-main} for terminal condition~$\varphi_\mu(T)=p$ and initial condition~$x_\mu(0)=0$. In the following, solving the linear system in~\eqref{equ:linear-system} and the optimality system in~\eqref{equ:optimality-system-odes} exactly will be referred to as solving the full-order model (\FOM{}). The main idea of the reduced order model introduced in the next section is to approximate the final time adjoint~$\varphi_\mu^*(T)$ by an element from a low-dimensional subspace of~$\setR^n$.

\section{Reduced order models for parametrized optimal control problems}\label{sec:reduced-order-models}
In this section we first introduce a reduced order model (ROM) based on a reduced basis approximation of the manifold~$\mathcal{M}\coloneqq\{\varphi_\mu^*(T):\mu\in\Params\}$ of optimal final time adjoint states. Afterwards, we show how to further accelerate the online phase by applying machine learning algorithms. We finally discuss an a posteriori error estimator for both reduced models.

\subsection{Reduced basis ROM}\label{subsec:rb-rom}
Assume that we are given a reduced basis~$\Phi^N=\{\varphi_1,\dots,\varphi_N\}$ for some~$\varphi_1,\dots,\varphi_N\in\setR^n$ and the respective reduced subspace~$X^N=\operatorname{span}(\Phi^N)\subset\setR^n$. In the adaptive model hierarchy described below, the reduced basis is built iteratively by starting with an empty basis and adding optimal final time adjoint states for certain parameters. A greedy procedure (see~\cite{devore2013greedy} for the theoretical background) to determine a reduced basis for the optimal control problem in~\eqref{equ:optimal-control-problem} was discussed in~\cite[Section~3]{kleikamp2023greedy}. Given the reduced space~$X^N$ and a parameter~$\mu\in\Params$, we compute the approximate final time adjoint~$\tilde{\varphi}_\mu^N\in X^N$ as
\[
	\tilde{\varphi}_\mu^N \coloneqq \argmin_{\tilde{\varphi}\in X^N}\norm{M\left(e^{A_\mu T}x_\mu^0-x_\mu^T\right) - \left(I+M\Gramian\right)\tilde{\varphi}}.
\]
In other words, we choose~$\tilde{\varphi}_\mu^N$ such that
\[
	\left(I+M\Gramian\right)\tilde{\varphi}_\mu^N = P_{Y_\mu^N}\Big((I+M\Gramian)\varphi_\mu^*(T)\Big),
\]
where~$P_{Y_\mu^N}$ denotes the orthogonal projection onto the space~$Y_\mu^N\coloneqq\left(I+M\Gramian\right)X^N$. This choice is motivated by the least squares solution of the linear system~\eqref{equ:linear-system} over the space~$X^N$.
\par
To compute the approximation~$\tilde{\varphi}_\mu^N\approx\varphi_\mu^*(T)$ in practice, one first computes the states~$x_i^\mu=(I+M\Gramian)\varphi_i$ for~$i=1,\dots,N$ (which essentially means to solve the optimality system in~\eqref{equ:optimality-system-main}). Afterwards, the matrix~$\bar{X}_\mu=[x_1^\mu\cdots x_N^\mu]\in\setR^{n\times N}$ can be assembled and the coefficients~$\alpha^\mu=(\alpha_1^\mu,\dots,\alpha_N^\mu)^\top\in\setR^N$ are derived as solutions of the linear system
\[
	\bar{X}_\mu^\top\bar{X}_\mu\alpha^\mu = \bar{X}_\mu^\top M\left(e^{A_\mu T}x_\mu^0-x_\mu^T\right).
\]
Having the coefficients at hand, the approximate final time adjoint is given as
\[
	\tilde{\varphi}_\mu^N = \sum\limits_{i=1}^N \alpha_i^\mu\varphi_i.
\]
The corresponding approximate optimal control~$\tilde{u}_\mu^N$ is computed by solving the optimality system in~\eqref{equ:optimality-system-main}. The reduced order model described in this section will be called reduced basis ROM (\RBROM{}) in the following.

\subsection{Machine learning ROM}\label{subsec:ml-rom}
The~\RBROM{} introduced in the previous subsection still involves several steps whose computational effort depends on the dimension~$n$ of the state space. In particular, computing~$x_i^\mu\in\setR^n$ for~$i=1,\dots,N$ corresponds to solving the adjoint and the state equation~$N$~times for every new parameter~$\mu\in\Params$. However, these computations are solely required to solve for the coefficients~$\alpha^\mu\in\setR^N$ with respect to the reduced basis. In the machine learning ROM (\MLROM{}) proposed in~\cite{kleikamp2023greedy}, instead of solving a linear system of equations for the coefficients, the map~$\pi_N\colon\Params\to\setR^N$, $\pi_N(\mu)\coloneqq\alpha^\mu$, from parameter to coefficients is approximated using machine learning algorithms. This idea is motivated by the approach first introduced in~\cite{hesthaven2018nonintrusive}. Given an approximation~$\hat{\pi}_N\colon\Params\to\setR^N$ of~$\pi_N$, the machine learning approximation~$\hat{\varphi}_\mu^N\in\setR^n$ of the optimal final time adjoint is defined as
\[
	\hat{\varphi}_\mu^N \coloneqq \sum\limits_{i=1}^N [\hat{\pi}_N(\mu)]_i\varphi_i.
\]
Similar to the~\RBROM{} in the previous section, the approximate optimal control~$\hat{u}_\mu^N$ is derived according to the optimality system from~\Cref{subsec:optimality-system}.
\par
The machine learning surrogate is trained in a supervised manner, i.e.~by means of training data consisting of parameters and corresponding coefficients. Several different machine learning algorithms are applicable in this scenario, the only requirement is that vector-valued functions can be approximated using training data (see~\cite[Section~4]{kleikamp2023greedy}). In our numerical experiment below we apply kernel methods (as discussed in more detail in~\cite[Section~4.3.2]{kleikamp2023greedy}) and in particular the \emph{vectorial kernel orthogonal greedy algorithm}~(VKOGA) as introduced in~\cite{santin2021kernel}.

\subsection{Residual based a posteriori error estimation}\label{subsec:error-estimation}
To estimate the error of the~\RBROM{} and the~\MLROM{} in an a posteriori manner, we consider the norm of the residual of the linear system in~\eqref{equ:linear-system}. To be more precise, given a parameter~$\mu\in\Params$ and an approximate final time adjoint~$p\in\setR^n$, the error estimate~$\eta_\mu(p)$ is defined as
\begin{align}\label{equ:error-estimator}
	\eta_\mu(p) \coloneqq \norm{M\left(e^{A_\mu T}x_\mu^0-x_\mu^T\right)-(I+M\Gramian)p}.
\end{align}
As proven in~\cite[Theorem~3.1]{kleikamp2023greedy}, it is possible to show that~$\eta_\mu$ is an efficient and reliable error estimator for the true error, i.e.~it holds
\begin{align}\label{equ:error-estimator-bounds}
	\norm{\varphi_\mu^*(T)-p} \quad \leq \quad \eta_\mu(p) \quad \leq \quad \norm{I+M\Gramian} \norm{\varphi_\mu^*(T)-p}.
\end{align}
This error estimator can be applied to both, the solution~$\tilde{\varphi}_\mu^N$ of the~\RBROM{} and the solution~$\hat{\varphi}_\mu^N$ of the~\MLROM{}. Evaluating the error estimator for some~$p\in\setR^n$ requires solving the adjoint equation once backwards in time with terminal condition~$p$, computing the corresponding control, and finally solving the state equation forward in time with zero initial condition for the state.

\section{Adaptive model hierarchy}\label{sec:adaptive-model-hierarchy}
In~\cite{haasdonk2023certified}, an adaptive model hierarchy for pa\-ra\-me\-tri\-zed PDEs was presented. The model hierarchy consists of a~\FOM{}, an~\RBROM{} and an~\MLROM{}. Both ROMs can be evaluated in terms of their a posteriori error using a residual-based error estimator. When the model hierarchy is queried for a new parameter~$\mu\in\Params$, first the~\MLROM{} is evaluated and its error compared to a prescribed error tolerance~$\varepsilon>0$. If the~\MLROM{}, which is the fastest of the three involved models, is sufficiently accurate, the machine learning approximation is returned. Otherwise, the~\RBROM{}, which takes more time to solve than the~\MLROM{} but is still faster than the~\FOM{}, is called. Also for the~\RBROM{}, the a posteriori error estimator is evaluated and the estimated error is compared to the tolerance~$\varepsilon$. If the~\RBROM{} is accurate enough, the reduced basis approximation is returned, if not, the solution provided by the~\FOM{} is computed. The~\FOM{} is typically much slower than both of the reduced models. The model hierarchy is constructed in such a way that, whenever possible, i.e.~their accuracy is sufficient, the faster to evaluate reduced order models are used and calls to the slower models are avoided.
\par
Furthermore, the model hierarchy is also adaptive in the sense that the reduced models are built while already querying the model hierarchy for different parameters. To be more precise, instead of using pre-trained reduced order models, one starts with an empty reduced basis such that the reduced order models can only return zero as solution. Whenever the~\FOM{} is called in the hierarchy, new training data for the~\RBROM{} is generated. Hence, calling the~\FOM{} improves the~\RBROM{}. Similarly, if the~\MLROM{} is not sufficiently accurate but the~\RBROM{} is, new training data for the~\MLROM{} is obtained by calling the~\RBROM{}. This way, the performance of both reduced models can be improved by solving the more expensive models in the hierarchy and thus benefit from more accurate solutions obtained by the slower models.

\section{Application of the model hierarchy to optimal control problems}\label{sec:application-model-hierarchy-optimal-control}
In the following we transfer the adaptive model hierarchy described in~\Cref{sec:adaptive-model-hierarchy} to parametrized optimal control problems by making use of the ROMs presented in~\Cref{subsec:rb-rom,subsec:ml-rom} together with the error estimator from~\Cref{subsec:error-estimation}. To this end, we assume that a desired accuracy~$\varepsilon>0$ of the computed (approximate) final time adjoint with respect to the optimal final time adjoint obtained by solving the~\FOM{} is prescribed.
\par
As a first step, the~\MLROM{} as introduced in~\Cref{subsec:ml-rom} is called to obtain the approximate final time adjoint~$\hat{\varphi}_\mu^N$. Afterwards, the error of the~\MLROM{} is estimated by evaluating~$\eta_\mu(\hat{\varphi}_\mu^N)$ where the error estimator~$\eta_\mu$ is defined in~\eqref{equ:error-estimator}. If the estimated error is smaller or equal to the tolerance~$\varepsilon$, the approximate optimal control~$\hat{u}_\mu^N$ is returned. If instead the~\MLROM{} was not accurate enough, the~\RBROM{} is solved for~$\tilde{\varphi}_\mu^N$. Similarly to the~\MLROM{}, the error estimate~$\eta_\mu(\tilde{\varphi}_\mu^N)$ is computed and compared to the prescribed tolerance~$\varepsilon$. If the~\RBROM{} is sufficiently accurate, i.e.~it holds~$\eta_\mu(\tilde{\varphi}_\mu^N)\leq\varepsilon$, the control~$\tilde{u}_\mu^N$ is returned. If even the~\RBROM{} is not accurate enough, the FOM is called and the optimal control~$u_\mu^*$ is returned.
\par
Due to the reliability of the error estimator~$\eta_\mu$ from~\eqref{equ:error-estimator}, i.e.~the first estimate in~\eqref{equ:error-estimator-bounds}, the result of the model hierarchy comes with a guaranteed accuracy. To be more precise, the error of the (approximate) final time adjoint used to compute the (approximate) optimal control that is returned by the model hierarchy is at most~$\varepsilon$ with respect to the optimal final time adjoint derived from the~\FOM{}.
\par
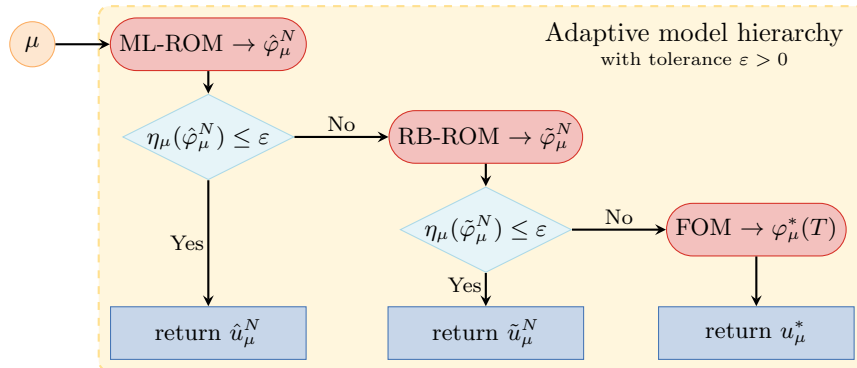
\begin{figure}[ht]
	\centering
	\begin{tikzpicture}[align=center, node distance=.3cm and 1.25cm]
		\node (mu) [input] {\hspace*{-2pt}$\mu$\hspace*{-2pt}};

		\node (mlm) [model, right=of mu, xshift=-15pt] {\MLROM{} $\rightarrow \hat{\varphi}_\mu^N$};
		\node (error-estimator-mlm) [error-estimator, below=of mlm] {$\eta_\mu(\hat{\varphi}_\mu^N)\leq\varepsilon$};
		\node (rb-rom) [model, right=of error-estimator-mlm] {\RBROM{} $\rightarrow \tilde{\varphi}_\mu^N$};
		\node (error-estimator-rb-rom) [error-estimator, below=of rb-rom] {$\eta_\mu(\tilde{\varphi}_\mu^N)\leq\varepsilon$};
		\node (fom) [model, right=of error-estimator-rb-rom] {FOM $\rightarrow \varphi_\mu^*(T)$};

		\node (output-mlm) at ($(error-estimator-mlm |- error-estimator-rb-rom.south)+(0,-.8cm)$) [output] {return $\hat{u}_\mu^N$};
		\node (output-rb-rom) at ($(error-estimator-rb-rom |- error-estimator-rb-rom.south)+(0,-.8cm)$) [output] {return $\tilde{u}_\mu^N$};
		\node (output-fom) at ($(fom |- error-estimator-rb-rom.south)+(0,-.8cm)$) [output] {return $u_\mu^*$};

		\draw[arrow] (mu) -- (mlm);
		\draw[arrow] (mlm) -- (error-estimator-mlm);
		\draw[arrow] (error-estimator-mlm) -- (output-mlm) node[pos=0.5, xshift=-8pt]{\footnotesize Yes};
		\draw[arrow] (error-estimator-mlm) -- (rb-rom) node[pos=0.5, yshift=5pt]{\footnotesize No};
		\draw[arrow] (rb-rom) -- (error-estimator-rb-rom);
		\draw[arrow] (error-estimator-rb-rom) -- (output-rb-rom) node[pos=0.35, xshift=-8pt]{\footnotesize Yes};
		\draw[arrow] (error-estimator-rb-rom) -- (fom) node[pos=0.5, yshift=5pt]{\footnotesize No};
		\draw[arrow] (fom) -- (output-fom);

		\background{output-mlm}{mlm}{output-fom}{output-fom}{bg-rb-rom};
		\node at (8.8,0) {Adaptive model hierarchy\\[-3pt]{\scriptsize with tolerance~$\varepsilon>0$}};
	\end{tikzpicture}
	\caption{Visualization of the adaptive model hierarchy applied to the parametrized optimal control setting.}
	\label{fig:adaptive-model-hierarchy}
\end{figure}

The way the adaptive model hierarchy determines a control when it is applied to parametrized optimal control problems and evaluated for a parameter~$\mu\in\Params$ is visualized in~\Cref{fig:adaptive-model-hierarchy}.

At this point, we recall again that every call to the~\RBROM{} in the model hierarchy generates new training data for the~\MLROM{}. Similarly, evaluating the FOM results in a new function for the reduced basis in the~\RBROM{}. Since extending the reduced basis results in a larger dimension of the reduced space, we also have to extend the previous~\MLROM{} which was created for a smaller reduced basis. One possibility to do so is to extend the previous training data by adding zeros for the new components as done in~\cite{haasdonk2023certified}. However, this strategy might result in an undesired bias towards zero for the~\MLROM{} training. To circumvent this issue, we follow a different policy in which we instead train a new machine learning surrogate for each individual coefficient in the reduced basis expansion. This means in particular that in the machine learning training no training data consisting of artificial zeros is used.
\par
We also emphasize that the main computational effort in the~\MLROM{} is spent for evaluating the error estimator. Computing the error estimate~$\eta_\mu(p)$ according to~\eqref{equ:error-estimator} for a given parameter~$\mu\in\Params$ and approximate final time adjoint~$p\in\setR^n$ still requires the solution of the adjoint equation and the primal state equation and therefore depends on the (large) dimension~$n$ of the state space. This computational effort in particular limits the efficiency of the~\MLROM{} when compared to the~\RBROM{} and will be discussed in more detail when investigating the numerical test case in the next section.

\section{Numerical experiment}\label{sec:numerical-experiment}
As a numerical example we consider a pa\-ra\-me\-trized heat equation where the parameter determines the heat conductivity of the underlying material as well as the target state. The two components of the control act on the Dirichlet boundary of the one-dimensional domain. The setting is similar to the one presented in~\cite[Section~6.2]{kleikamp2023greedy} and will be recalled briefly in the following (see~\cite{kleikamp2023greedy} for more details).
\par
Given a parameter~$\mu=[\mu_1,\mu_2]\in\Params\coloneqq[1,2]\times[0.5,1.5]\subset\setR^2$, the parametrized heat equation is given as
\begin{alignat*}{2}
	\partial_t v_\mu(t,y) - \mu_1\Delta v_\mu(t,y) &= 0 && \text{for }t\in[0,T],y\in\Omega, \\
	v_\mu(t,0) &= u_{\mu,1}(t) && \text{for }t\in[0,T], \\
	v_\mu(t,1) &= u_{\mu,2}(t) && \text{for }t\in[0,T], \\
	v_\mu(0,y) &= v_\mu^0(y) = \sin(\pi y) \qquad && \text{for }y\in\Omega,
\end{alignat*}
where~$v_\mu\colon[0,T]\times\Omega\to\setR$ with~$T=0.1$ denotes the state of the system. Furthermore, we denote by~$u_\mu(t)=\big[u_{\mu,1}(t),u_{\mu,2}(t)\big]^\top\in\setR^2$ for~$t\in[0,T]$ the control acting as Dirichlet boundary values on both ends of the one-dimensional domain~$\Omega=[0,1]$. In the optimal control problem we consider, the goal is to steer the system state at final time~$T$ close to the target state given by~$v_\mu^T(y)=\mu_2 y$ for~$y\in\Omega$. In this setting, the first component of the parameter determines the heat conductivity of the underlying material and the second component changes the slope of the target state. For the discretization of the system above, we use a second-order central finite difference scheme with~$200$ inner points in space and the Crank-Nicolson method with~$6000$ time steps for the time discretization.
\par
The adaptive model hierarchy is queried for~\numprint{10000} parameters from a uniformly distributed grid in the parameter set~$\Params$ that were randomly shuffled. The model hierarchy is applied to this system with a fixed tolerance of~$\varepsilon=10^{-4}$. Furthermore, the~\MLROM{} is trained whenever the reduced basis is extended or when five new training samples from the~\RBROM{} are collected. The experiment was performed on a dual socket compute server with two Intel(R) Xeon(R) Gold 6254 CPUs running at~3.10GHz and 36~cores in each CPU. The \texttt{Python} code for the experiment is available in~\cite{sourcecode} and can be used to reproduce the results shown below\footnote{The corresponding \texttt{GitHub}-repository containing the source code is available at~\url{https://github.com/HenKlei/ADAPTIVE-ML-OPT-CONTROL}}.
\par
\Cref{tab:heat-equation-results} summarizes the number of solves, number of error estimates, and timings for the~\FOM,~\RBROM{} and the~\MLROM{} within the adaptive model. We observe that the~\MLROM{} is sufficiently accurate in more than~99\% of the calls to the adaptive model hierarchy. The~\FOM{} was called~$4$~times in total which corresponds to a final reduced basis of size~$N=4$. We should remark at this point that the required amount of time for extending the reduced basis of the~\RBROM{} and for training the~\MLROM{} is negligibly small and was therefore omitted in the table. The advantage of using the adaptive model hierarchy instead of solely the~\FOM{} or the~\RBROM{} is reflected in the average time per solve of the three models (last column of~\Cref{tab:heat-equation-results}). The~\MLROM{} is about three times faster than the~\RBROM{} and about seventeen times faster than the~\FOM{}. This additional speedup of the~\MLROM{} pays off in view of the large number of evaluations of the model hierarchy. In particular, due to the error certification, the results obtained by the~\RBROM{} and the~\MLROM{} come with a guaranteed accuracy. However, we also observe that the speedup of the~\MLROM{} compared to the~\RBROM{} is only moderate. The reason for this observation is the relatively costly error estimation also for the~\MLROM{} as already discussed at the end of~\Cref{sec:application-model-hierarchy-optimal-control}. The speedup of the~\MLROM{} depends on the size of the reduced basis and becomes more pronounced for larger reduced bases.
\begin{table}[ht]
	\centering
		\begin{tabular}{l c c d{5.2} d{2.2}}
			\toprule
			{Model} & {\makecell{Number of\\solves}} & {\makecell{Number of\\error estimates}} & \mc{\makecell{Total time for error\\est.~and solving (s)}} & \mc{\makecell{Average time for error est.\\and solving per solve (s)}} \\ \midrule \midrule
			\FOM & \numprint{4} & {$-$} & 112.24 & 28.06 \\
			\RBROM & \numprint{65} & \numprint{69} & 299.26 & 4.60 \\
			\MLROM & \numprint{9931} & \numprint{10000} & 16{,}655.78 & 1.68 \\
			\bottomrule
	\end{tabular}
	\caption{Results of the numerical experiment for the heat equation using the adaptive model hierarchy.}
	\label{tab:heat-equation-results}
\end{table}
\par
\begin{figure}[ht]
	\centering
	\includegraphics{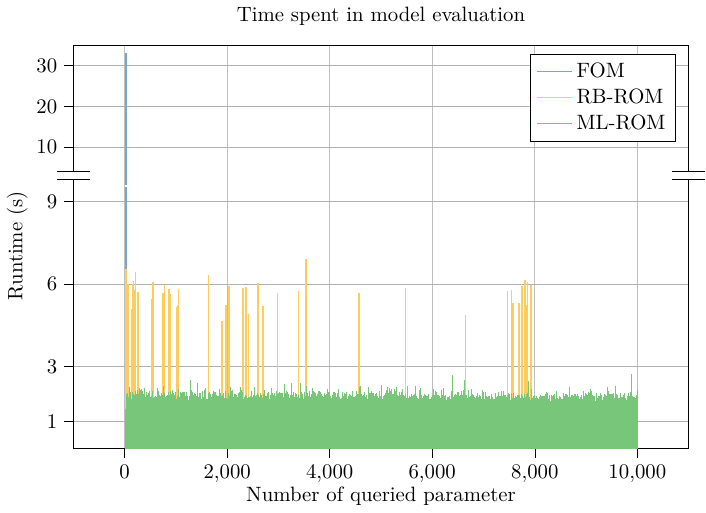}
	\caption{Performance of the adaptive model hierarchy in terms of the required times for error estimation and evaluation of the involved models when applied to a parametrized heat equation problem.}
	\label{fig:heat-equation-results}
\end{figure}

The calls to the models together with their runtimes over the queried parameters are shown in~\Cref{fig:heat-equation-results}. We observe that the~\FOM{} is called only for some of the first parameters and afterwards either the~\RBROM{} or the~\MLROM{} is sufficiently accurate and the~\FOM{} is never called again. In contrast, the model hierarchy falls back to the~\RBROM{} due to an insufficiently accurate~\MLROM{} for some parameters everywhere in the set of~\numprint{10000} parameters. Due to the a posteriori error estimation, such cases are detected and handled properly. The figure further shows that the~\MLROM{} is already used after a small number of parameters, i.e.~with a limited amount of training data an~\MLROM{} that is sufficiently accurate for several parameter values can be trained. In the last~\numprint{2000} calls to the model hierarchy, only the~\MLROM{} was used. As before, the runtimes for extending the reduced basis and training the~\MLROM{} are not shown in the plot since they are negligibly small.
\par
\begin{figure}[ht]
	\centering
	\includegraphics{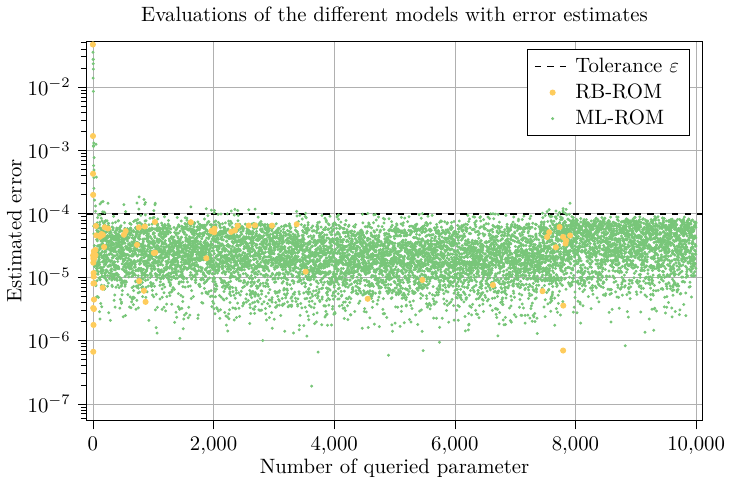}
	\caption{Error estimation of the~\RBROM{} and the~\MLROM{} in the adaptive model hierarchy when applied to a parametrized heat equation problem.}
	\label{fig:heat-equation-error-estimation-results}
\end{figure}

In~\Cref{fig:heat-equation-error-estimation-results} we further present the estimated errors of the~\RBROM{} and the~\MLROM{} over the queried set of parameters. As expected, for the first couple of parameters, the estimated error for both reduced models is above the desired tolerance~$\varepsilon$. After the four evaluations of the~\FOM{}, the~\RBROM{} is always sufficiently accurate. The estimated errors of the~\MLROM{} vary relatively strongly between about~$10^{-4}$ and~$10^{-6}$. For some parameters, the~\MLROM{} is not accurate enough and the~\RBROM{} is evaluated instead, but for several parameters the error of the~\MLROM{} is even about two orders of magnitude smaller than the prescribed tolerance.
\par
Altogether, applying the adaptive model hierarchy enables the certified approximate solution of the optimal control problem for~\numprint{10000} different parameters, where in the same time span only about~\numprint{608} calls to the~\FOM{} or~\numprint{3710} evaluations of the~\RBROM{} would have been possible.

\section{Conclusion and outlook}\label{sec:conclusion-outlook}
This work combines an adaptive model hierarchy introduced in~\cite{haasdonk2023certified} with reduced basis and machine learning ROMs for parametrized optimal control problems presented in~\cite{kleikamp2023greedy}. The main ingredient for the model hierarchy is the a posteriori error estimator available for both ROMs. By means of this error estimator, the model hierarchy is capable of providing certified results with guaranteed accuracy for every new parameter, while only using the fastest models whenever possible.
\par
The numerical example of a parametrized heat equation with controls acting on the Dirichlet boundary values shows how a combination of several layers of reduced order models can provide an additional speedup. At the same time, the error in the result of each query of the model hierarchy can be bounded by the prescribed tolerance due to the error certification. Hence, applying the model hierarchy results in a speedup while still maintaining the accuracy of the outputs.
\par
As future research directions, it might be of interest to investigate larger test cases and apply the model hierarchy in practical applications. Furthermore, one could investigate different choices of machine learning algorithms also in the setting of the adaptive model hierarchy. To further improve the performance of the model hierarchy, the evaluation of the error estimator should be sped up, since this constitutes the main remaining bottleneck in terms of computational efficiency of the overall procedure. The authors in~\cite{fabrini2018reduced} describe an approach based on a reduced basis ROM for the state and the adjoint equations which also accelerates the evaluation of the error estimator. However, using a reduction for the state and adjoint equations results in an additional error, which, to the best of our knowledge, has not been investigated so far.


\end{document}